\documentclass{amsart}
\usepackage{amsmath}
\usepackage{amsfonts}
\usepackage{amssymb,tikz}
\usepackage{tikz-cd}

\usepackage[all,  curve, frame, arc, color]{xy}
\usepackage{amscd}
\usepackage{mathrsfs}
\usepackage{soul}
\usepackage{enumitem}

\usepackage{hyperref}

\usepackage{array}  

\usepackage{proofread}
\usepackage[textsize=tiny]{todonotes}

\usepackage{sidenotes}
\usepackage{adjustbox}

\usepackage{color,soul}

\newcommand{\old}[1]{$^\star$#1} 
\newcommand{\olc}[1]{{#1}} 

\DeclareRobustCommand{\SkipTocEntry}[5]{}

\newtheorem{Thm}{Theorem}[section]

\theoremstyle{definition}

\newtheorem{Rem}[Thm]{Remark}

\newcommand{\Z}{\mathbb{Z}}

\newcommand{\Sc}{{\mathcal{S}}}
\newcommand{\Cc}{{\mathcal{C}}}

\newcommand{\A}{{\mathcal{A}}}

\usepackage{mathabx}

\makeatletter
\if@todonotes@disabled

\else

\fi
\makeatother

\begin{document}

\def\vg{\Lambda}	
\def\vd{\Omega}

\def\wl{\widehat{\lambda}}
\def\wg{\widehat{\omega}}

\def\T{\mathcal{T}}

\def\la{\langle}
\def\ra{\rangle}

\def\AL{\overline{\A}_{L}}
\def\AS{\overline{\A}_{S}}
\def\incl{\varphi}
\def\ideal{\mathcal{I}}

\def\G{\mathcal{G}}
\def\NN{\mathcal{N}}

\def\Stab{\operatorname{Stab}}
\def\supp{\operatorname{supp}}

\newcolumntype{C}{>{$}c<{$}}
\def\II{\mathcal{I}}
\def\JJ{\mathcal{J}}
\def\Gr{\operatorname{Gr}}
\def\tt{{\overline{t}}}

\title[Erratum: ``The integer cohomology algebra of toric arrangements'']{Erratum to \\ ``The integer cohomology algebra of toric arrangements''}
\author{Filippo Callegaro and Emanuele Delucchi}
\address[F. Callegaro]{
	Dipartimento di Matematica, University of Pisa, Italy.
}
\email{callegaro@dm.unipi.it}

\address[E. Delucchi]{University of Applied Arts and Sciences of Southern Switzerland.} 
\email{emanuele.delucchi@supsi.ch}

\begin{abstract}
We point out two errors in the paper ``The integer cohomology algebra of toric arrangements'', Adv.\ Math., Vol. 313, pp.~746--802, 2017. 

The main error concerns Theorem 4.2.17. In that theorem's proof, Diagram (8) does not commute in general but only under some restrictive hypotheses on the arrangement $\A$. This invalidates the description for the ring structure of $H^*(M(\mathcal{A});\mathbb{Z})$ given in Theorems A and B.  We refer to alternative descriptions of the cohomology ring $H^*(M(\mathcal{A});\mathbb{Z})$. 

The second error concerns Theorem 7.2.1. The claim does hold, but the proof is incorrect. We refer to a counterexample for the argument given in the proof and we provide references for a correct proof.
\end{abstract}
\maketitle
\vspace{-3em}
\setcounter{tocdepth}{1}
\tableofcontents

\noindent{\bf Note.} When referring to specific sections, theorems, lemmas, etc.\ in the original paper \cite{caldel17}, we will use a symbol \old{$\,$} in front of the number. Thus, for example, Theorem \old{3.4.7} means ``Theorem 3.4.7 in \cite{caldel17}''.

\section{Summary of corrections}
We start by an overview of the modifications in \cite{caldel17} that are required by the corrections in this Erratum. Some remarks for the corrections of \S \old{2--6} are given in Sections \ref{sec:algebra}. The details for the correction of \S \old{7} are given in Section \ref{sec:representations}.
\begin{itemize}[leftmargin=1.5em]
	\item[\S \old{2}] The rings $A(\A)$ and $B(\A)$ are isomorphic to a graded algebra associated to a filtration of $H^*(M(\A);\Z)$ induced by the Leray spectral sequence, but in general they are not isomorphic to the ring $H^*(M(\A);\Z)$ itself. The statements of Theorems~\old{A} and \old{B} should be replaced by Proposition 5.3.4, Theorem 5.3.5 and Theorem 5.3.7
	in \cite{erratum_long}, that hold when $\A$ is a real complexified
	toric arrangement. 
	\item[\S \old{3}] The results in this section hold without any change.
	\item[\S \old{4}] All statements until Equation~\old{(10)} in the proof of Theorem~\old{4.2.17} hold; the subsequent displayed equation may be false when $F_0\not \leq F_0'$. 
	
	\item[\S \old{5}]  Lemma~\old{5.1.2}, Theorem~\old{5.1.3} and Theorem~\old{5.1.5} hold. Corollary~\old{5.1.6} and the following statements are false.
	\item[\S \old{6}] Lemma~\old{6.1.2} and Theorem~\old{6.1.3} hold. Theorem~\old{6.2.4} is false. 
	We refer to \cite{a5,a5corrigendum} for a description of the cohomology of the complement in the general case.
	\item[\S \old{7}] The claim of Theorem~\old{7.2.1} holds, but the proof given in \cite{caldel17} is wrong. Moreover, Example \old{7.3} gives only a description of the graded ring associated to the filtration of the cohomology of the complement of the toric arrangement.
\end{itemize}

\section{Corrections of the cohomology algebra presentation}\label{sec:algebra}

The claim of Theorem~\old{4.2.17} in \cite{caldel17} does not hold for all choices of the basis chamber upon which relies  the construction of the subcomplexes $\Sc_L$. In particular, for some such choice the Diagram~\old{(8)} of \cite{caldel17} does not commute. This invalidates the description for the ring structure of $H^*(M(\A);\Z)$ given in Theorem~\old{A} and~\old{B}: in particular, the rings $A(\A)$ and $B(\A)$ are isomorphic to a graded algebra associated to a filtration of $H^*(M(\A);\Z)$ induced by the Leray spectral sequence, but in general they are not isomorphic to the ring $H^*(M(\A);\Z)$ itself. 

\begin{Rem}\label{rem:precisissimissimo}
The  flaw in the proof  of Theorem~\old{4.2.17}
lays in the computation leading from Equation~\old{(10)} to the next displayed equation. This computation does not hold in general, while it holds if $F_0\leq F_0'$. \footnote{Since Theorem~\old{4.2.17} assumes that $L \subseteq L'$, in principle one could choose faces $F_0,  F_0'$ such that $|F_0| = L, |F_0'| = L'$ and $F_0\leq F_0'$. However in  Definition~\old{4.2.16} we need to fix for every layer $L$ a face $F_0(L)$ of $\mathcal{A}_0$ such that $|F_0(L)| = L$. 
Hence it would be necessary to have a {\em global} choice of $F_0(L)$ such that for every $L,L'$ with $L \subseteq L'$ the relation  $F_0=F_0(L)\leq F_0(L') = F_0'$ holds.
In general such a global choice for $F_0(L)$ does not exist.}
\end{Rem}

\begin{Rem}\label{rem:quando_funziona}
The claim of Theorem~\old{4.2.17}  holds when it is possible to choose the facets $F_0$ that define the complexes $\mathcal{S}_L := \mathcal{S}_{F_0}$ in such a way that  $L=\vert F_0 \vert \subseteq \vert F'_0\vert=L'$  implies $F_0\leq F_0'$.  This is possible for instance when there is a fixed chamber $B_0$ of $\A_0$ such that for every layer $L$ the face $F_0$ associated to $L$ can be chosen as  $F_0=L_0\cap B_0$. A direct proof 
of the claim of Theorem~\old{4.2.17} in this case is given in \cite[Appendix A]{erratum_long}.
In particular, under the restrictive hypothesis on $\A$ that there exists a chamber $B_0$ of the arrangement $\A_0$ such that for every layer $L$ of $\A_0$ the support of the intersection $L \cap \overline{B}_0$ is $L$, the results of Theorem~\old{A} and~\old{B} 
hold.
\end{Rem}

In \cite{erratum_long} 
we show a workaround that allows  to provide a description of the cohomology ring $H^*(M(\A);\Z)$ as a subring of the direct sum $\oplus_{L \in \Cc} H^*(\Sc_{L};\Z)$ when $\A$ is a real complexified toric arrangement. Unfortunately this workaround 
does not apply to non real complexified toric arrangements.

For a general toric arrangement $\A$, not necessarily real complexified, a description of the ring $H^*(M(\A);\Z)$ in the style of the Orlik-Solomon algebra of hyperplane arrangements, obtained using other methods, can be found in \cite{a5,a5corrigendum}.

\section{Correction of the proof of  \cite[Theorem 7.2.1]{caldel17}}
\label{sec:representations}
In Section~\old{7} of \cite{caldel17} we investigate the dependency of our presentation of  from the combinatorial data. There, we claim the following result (where the last qualifier was implicit in the paragraphs preceding this theorem in \cite{caldel17}).
\begin{Thm}[{\cite[Theorem~\olc{7.2.1}]{caldel17}}]
If an arithmetic matroid with a basis of multiplicity $1$ is representable by a matrix $A$, then, if we fix such a basis, the matrix $A$ is unique up to sign reversal of the column vectors and up to a unimodular transformation from the left.
\end{Thm}


The proof given in the paper is not correct. As explained in \cite{Lenz}, the argument of the proof of case b) given in \cite{caldel17} fails for example for the matrix
$$
X = \left( \begin{array}{cccccc}
1 & 0 & 0 & 1 & 0 & 1\\
0 & 1 & 0 & 1 & 1 & 0\\
0 & 0 & 1 & 0 & 1  & -1
\end{array} \right)
$$
since it is not possible to make the bottom right entry positive preserving all other signs.

However, the claim of Theorem \old{7.2.1} of \cite{caldel17} is true: a correct proof follows from 
the results by Lenz \cite[Theorem~1.1]{Lenz} and 
in greater generality from Pagaria \cite[Theorem~3.5]{Pagaria}.

\bibliography{erratumbib}{}
\bibliographystyle{alpha}

\end{document}